\input amstex
\input amsppt.sty
\magnification=1200
\NoBlackBoxes
\def\interior{\operatorname{int}}
\def\PSH{\operatorname{PSH}}
\def\SH{\operatorname{SH}}
\def\diam{\operatorname{diam}}

\topmatter
\title
$L_h^2$-domains of holomorphy and the Bergman kernel
\endtitle
\author
Peter Pflug \& W\l odzimierz Zwonek
\endauthor
\abstract We give a characterization of $L_h^2$-domains of holomorphy with 
the help of the boundary behavior of the Bergman kernel and geometric
properties of the boundary, respectively.
\endabstract
\address
Carl von Ossietzky Universit\"at Oldenburg, Fachbereich Mathematik,
Postfach 2503, D-26111 Oldenburg, Germany
\endaddress
\address
Uniwersytet Jagiello\'nski, Instytut Matematyki, Reymonta 4, 30-059
Krak\'ow, Poland
\endaddress
\email pflug\@mathematik.uni-oldenburg, zwonek\@im.uj.edu.pl
\endemail
\thanks The second author was supported by the KBN grant No. 2 P03A 017 14.
While working on the paper the second author was at the stay at the Carl
von Ossietzky Universit\"at Oldenburg (Germany) supported by
DAAD.\endthanks \keywords Bergman kernel, $L_h^2$-domain of holomorphy,
(pluri)polar sets
\endkeywords
\endtopmatter

\document
For $\lambda_0\in\Bbb C$, $r>0$ we define
$\triangle(\lambda_0,r):=\{\lambda\in\Bbb C: |\lambda-\lambda_0|<r\}$. We
also put $E:=\triangle(0,1)$. Moreover, the set of all plurisubharmonic
(respectively, subharmonic) functions on an open set $D\subset\Bbb C^n$ is
denoted by $\PSH(D)$ (respectively, $\SH(D)$). We allow the
(pluri)subharmonic functions to be equal identically to $-\infty$ on
connected components of $D$.

Following \cite{Kli} for a domain $D\subset \Bbb C^n$ define
$$g_D(p,z):=\sup\{u(z)\},\quad p,z\in D,$$ where the supremum is taken
over all negative  $u\in\PSH(D)$ such that $u(\cdot)-\log||\cdot -p||$ is
bounded from above near $p$. We call the function $g_D(p,\cdot)$ {\it the
pluricomplex Green function } (with the logarithmic pole at $p$). We also
denote $$A_D(p;X):=\limsup_{\lambda\to 0}\frac{\exp(g_D(p,p+\lambda
z))}{|\lambda|},\quad p\in D,X\in\Bbb C^n.$$ Following \cite{Jar-Pfl} the
function $A_D$ is {\it the Azukawa pseudometric}.

For a boundary point $w$ of a bounded domain $D\subset\Bbb C$ we introduce
the notion of regularity. Namely, we say that $D$ is {\it regular } at $w$
if there exist a neighborhood $U$ of $w$ and a subharmonic function $u$
on $U\cap D$ with $u<0$ on $U\cap D$ and $\lim_{U\cap D\owns \lambda\to
w}u(\lambda)=0$.

A set $P\subset \Bbb C^n$ is called {\it pluripolar}
if for any point $z\in P$ there exist a connected neighborhood
$U=U(z)$ and a function $u\in\PSH(U)$, $u\not\equiv -\infty$,
such that $P\cap U\subset\{z\in U\:u(z)=-\infty\}$. In case
$n=1$ we call such a set $P$ {\it polar}. It is well known (cf. \cite{Kli},
Josefson theorem) that a set $P\subset\Bbb C^n$ is pluripolar if and only if 
there is a function $u\in PSH(\Bbb C^n)$, $u\not\equiv -\infty$, such that
$P\subset\{z\in\Bbb C^n:u(z)=-\infty\}$.

A bounded domain $D\subset \Bbb C^n$ is said to be {\it hyperconvex} if
there exists a negative and continuous plurisubharmonic exhaustion
function of $D$.

Denote the class of square integrable holomorphic functions on an open set
$D$ by $L_h^2(D)$. $L_h^2(D)$ is a Hilbert space with the standard scalar
product induced from $L^2(D)$. Let us recall the definition of {\it the
Bergman kernel}: $$
K_D(z):=\sup\{\frac{|f(z)|^2}{||f||_{L_h^2(D)}^2}:f\not\equiv 0,\;f\in
L_h^2(D)\}. $$ If $D$ is a bounded domain then $\log K_D$ is smooth and
strictly plurisubharmonic. Therefore, for a bounded domain $D$ one may
define {\it the Bergman metric $\beta_D$}: $$
\beta_D(z;X):=\sqrt{\sum_{j,k=1}^n\frac{\partial^2\log K_D(z)}{\partial
z_j\bar\partial z_k}X_j\bar X_k},\quad z\in D,X\in\Bbb C^n.$$ With the
help of the Bergman metric we obtain:
$$b_D(w,z):=\inf\{L_{\beta_D}(\alpha)\},\; w,z\in D$$ where
$L_{\beta_D}(\alpha)=\int_0^1\beta_D(\alpha(t);\alpha^{\prime}(t))dt$ and
the infimum is taken over all piecewise $C^1$-curves $[0,1]\mapsto D$. We
call $b_D$ the Bergman distance. If $(D,b_D)$ is a complete metric space
we say that $D$ is {\it Bergman complete}.

A domain $D\subset\Bbb C^n$ is called {\it a (an $L_h^2$-)domain of
holomorphy} if there are no domains $D_0,D_1\subset\Bbb C^n$ with
$\emptyset\neq D_0\subset D_1\cap D$, $D_1\not\subset D$ such that for any
$f\in\Cal O(D)$ ($f\in L_h^2(D)$) there exists an $\tilde f\in\Cal O(D_1)$
with $\tilde f=f$ on $D_0$.

Let us recall several results concerning the above mentioned notions,
which show a close relationship between the theory of square integrable
holomorphic functions and that of the pluripotential theory.

For a bounded pseudoconvex domain $D$ consider the following properties:
$$ \gather
  \text{$D$ is hyperconvex},\tag{1}\\
  \text{for any $w\in\partial D$, $\lim_{D\owns z\to w}K_D(z)=\infty$,}
  \tag{2}\\
  \text{$D$ is Bergman complete},\tag{3}\\
  \text{$D$ is an $L_h^2$-domain of holomorphy}.\tag{4}
\endgather
$$ All the relations between the properties \thetag{1}--\thetag{4} are
known. Namely, (\thetag{1}$\implies$\thetag{2}) (see \cite{Ohs~1}),
(\thetag{1}$\implies$\thetag{3}) (see \cite{B\l o-Pfl}, \cite{Her}),
(\thetag{3}$\implies$\thetag{4}). The implication
(\thetag{2}$\implies$\thetag{1}) does not hold in general (take the
Hartogs triangle in $\Bbb C^2$ or consider some one dimensional
Zalcmann-type domains -- see \cite{Ohs~1}). The one-dimensional
counterexample to the implication (\thetag{3}$\implies$\thetag{1}) is
given in \cite{Chen~1}. Recall that any bounded pseudoconvex fat domain is an
$L_h^2$-domain of holomorphy (see \cite{Pfl}). Thus the Hartogs triangle
is an $L_h^2$-domain of holomorphy in $\Bbb C^2$ which is not
Bergman complete. Moreover, there exists also a fat domain in the 
complex plane that is not Bergman complete (see \cite{Jar-Pfl-Zwo}). Thus, the
implication (\thetag{4}$\implies$\thetag{3}) does not hold even for fat 
pseudoconvex domains. In dimension one the
implication (\thetag{2}$\implies$\thetag{3}) holds (see \cite{Chen~2}) but
in higher dimension this is no longer the case (take the Hartogs triangle
once more). As far as the implication (\thetag{3}$\implies$\thetag{2}) is
concerned one may find a counterexample already in dimension one (see
\cite{Zwo~2}).

Let us have a closer look at the last example. The domains being
counterexamples belong to the following class of domains: $$
D:=E\setminus(\bigcup_{j=1}^{\infty}\bar\triangle(z_j,r_j)\cup\{0\}),$$
where $z_j\to 0$, $r_j>0$, $\bar\triangle(z_j,r_j)\subset
E\setminus\{0\}$,
$\bar\triangle(z_j,r_j)\cap\bar\triangle(z_k,r_k)=\emptyset$, $j\neq k$.
It is easy to see that for any $w\in\partial D$, $w\neq 0$ we have
$\lim_{D\owns z\to w}K_D(z)=\infty$. The point is that the sequences can
be chosen so that $\liminf_{D\owns z\to 0}K_D(z)<\infty$ and the domain is
still Bergman complete. On the other hand one may easily see that
$\limsup_{z\to 0}K_D(z)=\infty$. So the natural problem arises whether one
may construct an example of a Bergman complete domain such that for some
$w\in\partial D$ we have $\limsup_{z\to w}K_D(z)<\infty$. Below we show
that this is impossible. Let us write down explicitly the condition we are
interested in (as some kind of appendix to properties
(\thetag{1}--\thetag{4})): $$\text{for any $w\in\partial D$ we have
$\limsup_{D\owns z\to w}K_D(z)=\infty$.}\tag{5}$$

The main aim of this paper is to present the following characterizations
of $L_h^2$-domains of holomorphy.

\proclaim{Theorem 1} Let $D$ be a bounded pseudoconvex domain in $\Bbb
C^n$. Then \thetag{4} is equivalent to \thetag{5}.
\endproclaim

Making use of Theorem 1 and a result of A. Sadullaev we also get the
following characterization of bounded $L_h^2$-domains of holomorphy.

\proclaim{Theorem 2} Let $D$ be a bounded pseudoconvex domain. Then $D$ is
an $L_h^2$-domain of holomorphy if and only if for any $w\in\partial D$
and for any neighborhood $U$ of $w$ the set $U\setminus D$ is not
pluripolar.
\endproclaim

Before proving Theorem 1 let us recall some properties of the notions that
we have just defined and that we need in the sequel.

First we list a number of properties of polar sets that we shall use (see
\cite{Ran}, \cite{Con}).

Let $D$ be an open set in $\Bbb C$ and let $K\subset D$ be a polar set
relatively closed in $D$. Then

-- if $D$ is additionally connected then so is $D\setminus K$,

-- for any $\lambda\in D$ and for any $0<s$ with
$\triangle(\lambda,s)\subset\subset D$ there is an $s<r$ with
$\triangle(\lambda,r)\subset\subset D$ and
$\partial\triangle(\lambda,r)\cap K=\emptyset$,

-- for any $f\in L_h^2(D\setminus K)$ there is an $\tilde f\in O(D)$ such
that $\tilde f_{|D\setminus K}=f$.

There is also a precise description of $L^2_h$-domains of holomorphy in
$\Bbb C$. 

\proclaim{Theorem 3 {\rm (see \cite{Con}, Theorem 9.9, p. 351)}}
Let $D$ be a bounded domain in $\Bbb C$ and let $z\in\partial D$. Then
there is an open neighborhood $U$ of $z$ such that any $f\in L_h^2(D)$
extends holomorphically to $D\cup U$ if and only if there is a
neighborhood $V$ of $z$ such that the set $V\setminus D$ is polar.
\endproclaim

One may easily get from Theorem 3 the following description of
$L_h^2$-domains of holomorphy in $\Bbb C$. 

\proclaim{Theorem 4} Let $D$ be a bounded domain in $\Bbb C$. Then $D$ 
is an $L_h^2$-domain of holomorphy
iff for any $w\in\partial D$ and for any neighborhood $U$ of $w$ the set
$U\setminus D$ is not polar.
\endproclaim

Note that Theorem 2 is the exact more dimensional counterpart of Theorem
4.

Let us recall now some basic properties of regular points and the Green
function. For a domain $D\subset\Bbb C^n$ we have
$g_D(p,\cdot)\in\PSH(D)$, $g_D(p,\cdot)<0$. A bounded domain $D$ is
hyperconvex iff $g_D(p,\cdot)$ is a continuous exhaustive function of $D$.

In the case of bounded planar domains it is well-known that the Green
function is symmetric (as the function of two variables) and
$g_D(p,\cdot)$ is harmonic on $D\setminus\{p\}$. Moreover, a point
$w\in\partial D$ is regular iff for some (any) $p\in D$ $g_D(p,\lambda)\to
0$ as $D\owns\lambda\to w$. Consequently, a bounded domain $D\subset\Bbb
C$ is hyperconvex iff any point from its boundary is regular. The set of
irregular points of any bounded domain in $\Bbb C$ is polar.

Below we shall need some estimate for the Bergman kernel in the
one-dimensional case that will enable us to prove Theorem 1 in dimension
one.

\proclaim{Theorem 5 {\rm (see \cite{Ohs~2})}} Let $D$ be a domain in $\Bbb
C$. Then there is a positive constant $C$ such that $$\sqrt{K_D(z)}\geq
CA_D(z;1),\;z\in D.$$
\endproclaim

Our first aim is to obtain the following exhaustion property of the
Bergman kernel at regular points.

\proclaim{Proposition 6} Let $D$ be a bounded domain in $\Bbb C$. Assume
that $w\in\partial D$ is a regular point. Then $$K_D(z)\to\infty$$ as
$D\owns z\to w$.
\endproclaim
\demo{Proof} In view of Theorem 5 it is sufficient to show that $$r(p)\to
0 \text{ as $p\to w$},\tag{6}$$ where $r:=r(p):=\diam D(p)$, $D(p):=\{z\in
D:g_D(p,z)<-1\}$. In fact, assuming the last property we get (see
\cite{Zwo~1}) $$A_D(p;1)=eA_{D(p)}(p;1)\geq
eA_{\triangle(p,r)}(p;1)=\frac{e}{r}\to\infty \text{ as $p\to w$}.$$
Suppose that \thetag{6} does not hold. Then one easily finds an
$\epsilon>0$, sequences $D\owns p_{\nu}\to w$ and $D\owns z_{\nu}\to
z\in\bar D$ such that $|p_{\nu}-z_{\nu}|\geq \epsilon$ and
$g_D(p_{\nu},z_{\nu})<-1$. Taking $\tilde D:=D\cup V$, where $V$ is some
small disc around $z$ such that $w\not\in\bar V$, we get
$g_D(p_{\nu},z_{\nu})\geq g_{\tilde D}(p_{\nu},z_{\nu})$ and $z\in\tilde
D$. In other words, it is sufficient to show that $g_{\tilde
D}(p_{\nu},z_{\nu})\to 0$. But because of the pointwise convergence of
$g_{\tilde D}(p_{\nu},\cdot)=g_{\tilde D}(\cdot,p_{\nu})$ to $0$ (as
$\nu\to\infty$), the harmonicity of $g_{\tilde D}(p_{\nu},\cdot)$ near $z$
and the Vitali theorem, we conclude that $g_{\tilde D}(p_{\nu},\cdot)$ tends
uniformly to $0$ on some neighborhood of $z$, which finishes the proof.
\qed
\enddemo

\subheading{Remark 7} In view of property \thetag{6} it follows from the
estimates in \cite{Die-Her} that for any bounded domain in $\Bbb C$ the
convergence $\beta_D(z;1)\to\infty$ as $z\to w\in\partial D$ holds for any
regular point $w\in\partial D$.

\proclaim{Lemma 8} Let $D$ be a bounded domain in $\Bbb C$, $w\in\partial
D$. Then the following conditions are equivalent: $$ \gather
\limsup_{D\owns z\to w}K_D(z)<\infty,\tag{7}\\ \text{there is an open
neighborhood $U$ of $w$ such that the set $U\setminus D$ is
polar.}\tag{8}
\endgather
$$
\endproclaim
\demo{Proof} Let us make at first a general remark. Namely, the condition
$U\setminus D$ is polar is equivalent to the condition $U\cap\partial D$
is polar.

(\thetag{8}$\implies$\thetag{7}). If $U$ satisfies \thetag{8} then without
loss of generality one may assume that $K:=U\cap\partial D\subset\subset
U$. So there is a domain $\tilde D$ with $D=\tilde D\setminus K$,
$w\in\tilde D$, where $K$ is a compact polar set. Then $$
L_h^2(D)=L_h^2(\tilde D)_{|D} $$ and, consequently, $$ K_D={K_{\tilde
D}}_{|D}, $$ which implies \thetag{7}.

(\thetag{7}$\implies$\thetag{8}). Suppose that for any neighborhood $U$
of $w$ the set $U\cap\partial D$ is not polar. Then there is a sequence
$w_{\nu}\to w$, $w_{\nu}\in\partial D$, such that $D$ is regular at
$w_{\nu}$. In view of Proposition 6 we have $K_D(z)\to\infty$ as $D\owns
z\to w_{\nu}$, which easily finishes the proof. \qed
\enddemo

\proclaim{Lemma 9} Let $D$ be a domain in $\Bbb C^n$, $n\geq 2$. Fix
$0<r<t$. For any $z^{\prime}\in \Bbb C^{n-1}$ define 
$A(z^{\prime}):= \{z_n\in tE:(z^{\prime},z_n)\in D\}
:=tE\setminus K(z^{\prime})$. Assume that
$K(0^{\prime})$ is polar and there is a neighborhood $0^{\prime}\in V$
such that for almost any $z^{\prime}\in V$ (with respect to the
$(2n-2)$-dimensional Lebesgue measure) the set $K(z^{\prime})$ is polar.
Then there is a neighborhood $0^{\prime}\in V^{\prime}\subset V$ such that
for any $f\in L_h^2(D)$ there exists a function $F\in\Cal
O(V^{\prime}\times rE)$ with $F=f$ on $(V^{\prime}\times rE)\cap D$.
\endproclaim
\demo{Proof} Because $K(0^{\prime})$ is polar there is an $s$ with
$0<r<s<t$ such that $K(0^{\prime})\cap
\partial(sE)=\emptyset$. Then there is a neighborhood
$0^{\prime}\in V^{\prime}\subset V$ such that for any $\zeta^{\prime}\in
V^{\prime}$ we have $K(\zeta^{\prime})\cap\partial(sE)=\emptyset$.

Define $$F(\zeta^{\prime},z_n):=\frac{1}{2\pi i}\int_{\partial(sE)}
\frac{f(\zeta^{\prime},\lambda)d\lambda}{\lambda-z_n},\quad
(\zeta^{\prime},z_n)\in V^{\prime}\times sE. $$ Then $F$ is a holomorphic
function on $V^{\prime}\times sE$.

On the other hand because of the square integrability of $f$, the Fubini
Theorem and because of the assumptions of the lemma, for almost all
$\zeta^{\prime}\in V^{\prime}$ (with respect to the $(2n-2)$-dimensional
Lebesgue measure) the function $f(\zeta^{\prime},\cdot)\in
L_h^2(tE\setminus K(\zeta^{\prime}))$ and $K(z^{\prime})$ is polar. Since
closed polar sets are removable for $L_h^2$-functions, for almost all
$\zeta^{\prime}\in V^{\prime}$ the function $f(\zeta^{\prime},\cdot)$
extends to a holomorphic function on $tE$. So the Cauchy formula applies
and we obtain the equality $f(\zeta^{\prime},z_n)=F(\zeta^{\prime},z_n)$,
$(\zeta^{\prime},z_n)\in (V^{\prime}\times sE)\cap D$ for almost all
$\zeta^{\prime}\in V^{\prime}$. Since equality holds on a dense subset of
$(V^{\prime}\times sE)\cap D$, the equality holds on the whole set, which
finishes the proof. \qed
\enddemo

Before we start the proof of Theorem 1 let us formulate, in the form that
we need, the most powerful tool we shall use, namely the Ohsawa-Takegoshi
extension theorem.

\proclaim{Theorem 10 {\rm(see \cite{Ohs-Tak})}} Let $D$ be a bounded
pseudoconvex domain in $\Bbb C^n$ and let $L$ be a complex line. Then
there is a constant $C>0$ such that for any $f\in L_h^2(D\cap L)$ there is
an $F\in L_h^2(D)$ with $||F||_{L_h^2(D)}\leq C||f||_{L_h^2(D\cap L)}$ and
$F_{|D\cap L}=f$.
\endproclaim

Note that Theorem 10 directly leads to the following  inequality for the
Bergman kernel: $$K_{D\cap L}(z)\leq C^2K_D(z),\;z\in D\cap L.$$ This
inequality will be often used in our next considerations. Note only that
the set $D\cap L$ on the left-hand side of the above inequality is open
(as a subset of $\Bbb C$) but not necessarily connected.

We now prove our main result. 

\demo{Proof of Theorem 1} First note that
the result for $n=1$ follows from Theorem 4 and Lemma 8, so assume that
$n\geq 2$.

(\thetag{5}$\implies$\thetag{4}). Suppose that $D$ is not an
$L_h^2$-domain of holomorphy. Then there are a polydisc $P\subset D$ with
$\partial P\cap
\partial D\neq\emptyset$ and a polydisc $P\subset\subset\tilde P$ such
that for every function $f\in L_h^2(D)$ there is a function $\hat f\in
H^{\infty}(\tilde P)$ with the property $f=\hat f$ on $P$.

We claim that for any $z\in P$ and for any complex line $L$ passing
through $z$ we have $$L\cap D\cap \tilde P=(L\cap \tilde P)\setminus K(z),
\text{ where $K(z)$ is a polar set.}$$

Suppose that $L\cap D\cap \tilde P=(L\cap \tilde P)\setminus K(z)$, where
$K(z)$ is not a polar set. Then choose a compact non-polar set
$K^{\prime}\subset K(z)\subset (L\cap\tilde P)\setminus D$ such that
$V_0=L\setminus \hat K^{\prime}$ ($\hat K^{\prime}$ denotes the polynomial
hull of $K^{\prime}$) contains $L\cap P$. Then there is a function $f\in
L_h^2(V_0)$ which does not extend holomorphically through $\hat
K^{\prime}$ (cf. Theorem 3). Let $\{V_j\}_{j=1}^N$, where $0\leq N\leq
\infty$ be the family of bounded components of $L\setminus K^{\prime}$.
We, additionally, let $f$ be identically $0$ on $\bigcup_{j=1}^NV_j$.

In view of the Ohsawa-Takegoshi extension theorem there exists an $F\in
L_h^2(D)$ such that $F_{|L\cap D}=f_{|L\cap D}$. But then there is a
function $\hat F\in H^{\infty}(\tilde P)$ such that $\hat F_{|P}=F_{|P}$.
Consequently, $\hat F_{|L\cap\tilde P}$ is a holomorphic extension of
$f_{|L\setminus \hat K^{\prime}}$ through $\hat K^{\prime}$ --
contradiction.

It follows from the above claim that $\tilde P\cap D$ is connected.
Consequently,
for any function $f\in L_h^2(D)$ its (uniquely determined) extension
$\hat f\in H^{\infty}(\tilde P)$ satisfies the equality $f=\hat f$ on
$D\cap\tilde P$.

Consider the following normed space $$ A:=\{(f,\hat f):f\in
L_h^2(D)\}\subset L_h^2(D)\times H^{\infty}(\tilde P) $$ with the norm
$||(f,\hat f)||:=||f||_{L_h^2(D)}+||\hat f||_{H^{\infty} (\tilde P)}$. It
is easily seen that $A$ is a Banach space. Consider the following mapping:
$$ \pi: A\owns (f,\hat f)\mapsto f\in L_h^2(D). $$ Then $\pi$ is a
one-to-one surjective continuous linear mapping. Hence, in view of the
Banach open mapping theorem, $\pi^{-1}$ is a continuous linear mapping. In
other words, there is a constant $C>1$ such that $$ ||(f,\hat f)||\leq
C||f||_{L_h^2(D)},\quad f\in L_h^2(D);$$ in particular, $||\hat
f||_{H^{\infty}(\tilde P)}\leq C||f||_{L_h^2(D)}$. Consequently, $$
\sup_{z\in\tilde P\cap D}K_D(z)=\sup\{\frac{|f(z)|^2}{||f||_{L_h^2(D)}^2}:
z\in\tilde P\cap D,\,f\not\equiv 0,\;f\in L_h^2(D)\}\leq C^2, $$ which
contradicts \thetag{5} for any $w\in\partial P\cap\partial
D\neq\emptyset$.

(\thetag{4}$\implies$\thetag{5}). Fix $w\in\partial D$.

Let us first consider the case $w\not\in\interior(\bar D)$. Then there is
a sequence $z_{\nu}\to w$, $z_{\nu}\not\in\bar D$. Let $B_{\nu}$ be the
largest open ball centered at $z_{\nu}$ disjoint from $\bar D$. Choose
$w_{\nu}\in\partial B_{\nu}\cap\partial D$. Obviously, $w_{\nu}\to w$.
Note that for any $\nu$, $D$ satisfies at $w_{\nu}$ 'the outer cone
condition' (see \cite{Pfl}). Therefore, for any $\nu$ we have
$\lim_{D\owns z\to w_{\nu}}K_D(z)=\infty$ (see \cite{Pfl}), which easily
implies $\thetag{5}$.

Assume now that $w\in\interior(\bar D)$. Suppose that \thetag{5} does not
hold at $w$. Then there is a polydisc $P$ with centre at $w$ such that
$\sup\{K_D(z):z\in D\cap P\}<\infty$. Without loss of generality we may
assume that $P\subset\subset \interior(\bar D)$. Consider any complex line
$L$ intersecting $P$. We claim that $L\cap P\cap D$ is equal to $(L\cap
P)\setminus K$, where $K$ is a polar set or $K=L\cap P$. In fact if this
were not the case then $\sup_{z\in L\cap P\cap D}K_{L\cap D}(z)= \infty$
(the Bergman kernel is here understood as that of a one-dimensional set)
(use Lemma 8) and, consequently, in view of the Ohsawa-Takegoshi extension
Theorem we get $\sup_{z\in L\cap P\cap D}K_D(z)=\infty$ -- contradiction.

Note that there is a complex line $L$ passing through $w$ such that $L\cap
P\cap D$ is not empty. Assume that $w=0$. Making a linear change of
coordinates and shrinking $P$ if necessary we may assume that $P=E^n$ and
that $\{\lambda\in E: (0,\ldots,0,\lambda)\in D\}$ is not empty.

Therefore, the assumptions of Lemma 9 are satisfied (with some
neighborhood $V\subset E^{n-1}$ of $0^{\prime}\in\Bbb C^{n-1}$) and there is a
neighborhood $0^{\prime}\in V^{\prime}\subset E^{n-1}$ such that for any $f\in
L_h^2(D)$ there is a function $F\in\Cal O(V^{\prime}\times\frac{1}{2}E)$
with $F=f$ on $(V^{\prime}\times\frac{1}{2}E)\cap D$ -- contradiction.
\qed \enddemo

\demo{Proof of Theorem 2} Because of Theorem 4 we may assume
that $n\geq 2$.

($\Rightarrow$). Suppose that for some
$w\in\partial D$ there is a polydisc $P$ such that $P\setminus D$ is
pluripolar. Let $u\in\PSH(P)$ be such that $u\not\equiv-\infty$ and
$P\setminus D\subset\{u=-\infty\}$.  Take a nonempty open set $U\subset
D\cap P$ and consider all complex lines connecting $w$ with some point
from $U$. It is easy to see that there is a complex line $L$ such that
$u\not\equiv-\infty$ on $L\cap P$. Assume that $w=0$. Making linear change
of coordinates and shrinking $P$, if necessary, we may assume that $P=E^n$
and $\{z_n\in E:(0^{\prime},z_n)\not\in D\}$ is polar. Because of the
local integrability of $u$, for almost any $z^{\prime}\in E^{n-1}$ (with
respect to the $(2n-2)$-dimensional Lebesgue measure) the function
$u(z^{\prime},\cdot)$ is not identically equal to $-\infty$ on $E$.
Consequently, for almost any $z^{\prime}\in E^{n-1}$ the set $\{z_n\in
E:(z^{\prime},z_n)\not\in D\}$ is polar. Applying Lemma 9 we obtain the
existence of an open set $0\in Q$ such that for any $f\in L_h^2(D)$ there
exists an $F\in\Cal O(Q)$ with $f=F$ on $D\cap Q$ -- contradiction.

($\Leftarrow$). Suppose that the implication does not hold, so in view of
Theorem 1 there is a $w\in\partial D$ such that $\limsup_{D\owns z\to
w}K_D(z)<\infty$. In other words there is a polydisc $P$ with centre at
$w$ such that $\sup_{z\in D\cap P}K_D(z)<\infty$.

First note that for any complex line $L$ with $L\cap P\neq\emptyset$ we
have $L\cap P\cap D=\emptyset$ or $L\cap P\cap D=(L\cap P)\setminus K$,
where $K$ is a polar set. Actually, if there were $L$ such that $L\cap
P\cap D=(L\cap P)\setminus K$, where $K\neq L\cap P$ and $K$ is not polar,
then for some $U\subset\subset L\cap P$, $\sup_{z\in U\cap D}K_{D\cap
L}(z)=\infty$ (use Lemma 8). Therefore, in view of the Ohsawa-Takegoshi
theorem, $\sup_{z\in U\cap D}K_D(z)=\infty$ -- contradiction.

Consequently, one may apply a result of A. Sadullaev (see \cite{Sad~2} and
also \cite{Sad~1}) to get that the set $P\setminus D$ is pluripolar --
contradiction.
\qed
\enddemo

It follows from the reasoning in proofs of Theorems 1 and 2 that the
following more dimensional counterpart of Lemma 8 holds 

\proclaim{Lemma 11} Let $D$ be a bounded pseudoconvex domain and let 
$w\in\partial D$.
Then $\limsup_{D\owns z\to w}K_D(z)<\infty$ if and only if for any
neighborhood $U$ of $w$ the set $U\setminus D$ is pluripolar.
\endproclaim

The known examples of $L_h^2$-domains of holomorphy include, among others,
bound\-ed pseudoconvex fat domains and bounded pseudoconvex balanced
domains. The characterization of $L_h^2$-domains of holomorphy given by us
yields many examples of such domains. Below we give one example of a new
class of domains having this property.

For a bounded pseudoconvex domain $D\subset\Bbb C^n$ we define the
following Hartogs domain with $m$-dimensional balanced fibers:
$$G_D:=\{(w,z)\in\Bbb C^{n+m}:H(z,w)<1\},$$ where $\log H$ is
plurisubharmonic on $D\times\Bbb C^m$, $H(z,\lambda w)=|\lambda|H(z,w)$,
$(z,w)\in D\times\Bbb C^m$, $\lambda\in\Bbb C$, and $G_D$ is bounded (i.e.
$H(z,w)\geq C||w||$ for some $C>0$, $(z,w)\in D\times\Bbb C^m$). In such a
situation $G_D$ is a bounded pseudoconvex domain.

\proclaim{Proposition 12} Let $D$ be a bounded $L_h^2$-domain of
holomorphy. Then $G_D$ (with notation as above) is an $L_h^2$-domain of
holomorphy.\endproclaim \demo{Proof} Let us take $(z^0,w^0)\in\partial
G_D$. If $z^0\in D$ then
$\lim_{G_D\owns(z,w)\to(z^0,w^0)}K_{G_D}(z,w)=\infty$ (use Theorem 3.1(i)
from \cite{Jar-Pfl-Zwo}).

Assume now that $z^0\in\partial D$. Let $V$ be any neighborhood of
$(z^0,w^0)$. In view of Lemma 11 and Theorem 1 it is sufficient to show
that $V\setminus G_D$ is not pluripolar.  Without loss of generality we
may assume that $V=V_1\times V_2\subset\Bbb C^{n+m}$. Because $D$ is an
$L_h^2$-domain of holomorphy Theorem 2 applies and $V_1\setminus D$ is not
pluripolar. Since $V\setminus G_D\supset (V_1\setminus D)\times V_2$ and
the latter set is not pluripolar, the proof is finished.
\qed
\enddemo

\subheading{Acknowledgment} The authors would like to thank A. Edigarian
for drawing their attention to a result of A. Sadullaev.

\enddocument